\documentclass[10pt]{amsart}
\usepackage{amsmath}
\usepackage{amssymb}
\usepackage{amsthm}
\usepackage{psfig}

\newcommand{\ignore}[1]{\relax}

\newcommand{\C}{\mathbb C}
\newcommand{\R}{\mathbb R}
\newcommand{\Z}{\mathbb Z}
\newcommand{\N}{\mathbb N}

\newcommand{\U}{\mathcal U}

\newcommand{\HH}{\mathcal H}
\newcommand{\LLL}{\mathcal L}

\newtheorem{thm}{Theorem}

\theoremstyle{definition}

\newtheorem{exa}{Example}
\newtheorem{exam}{Special Case}

\newtheorem{rmk}{Remark}

\newcommand{\tor}{(\C^*)^{2}}
\newcommand{\rtor}{(\R^*)^{2}}

\newcommand{\dd}{\partial}

\newcommand{\cp}{{\mathbb C}{\mathbb P}}
\newcommand{\rp}{{\mathbb R}{\mathbb P}}
\newcommand{\pp}{\mathbb P}

\newcommand{\Int}{\operatorname{Int}}
\renewcommand{\setminus}{\smallsetminus}

\begin{document}

\title
%[Curves in toric surfaces]
{Counting curves via lattice paths in polygons}
\author{Grigory Mikhalkin}
\thanks{The author is partially supported by the NSF}
\address{Dept of Mathematics\\
Univ of Utah\\
Salt Lake City, UT 84112, USA}
\address{St. Petersburg Branch of Steklov Mathematical Institute,
Fontanka 27, St. Petersburg, 191011 Russia}
\email{mikhalkin@math.utah.edu}
%\maketitle

\begin{abstract}
This note presents a formula for the enumerative invariants
of arbitrary genus in toric surfaces.
The formula computes the number of curves of a given genus through
a collection of generic points in the surface.
%The number of points
%is chosen so that the number of relevant curves is finite.
The answer is given
%as a number (with multiplicities)
in terms of certain lattice paths in the relevant Newton polygon.
If the toric surface is ${\mathbb P}^2$ or ${\mathbb P}^1\times
{\mathbb P}^1$ then the invariants under consideration
coincide with the Gromov-Witten invariants.
The formula gives a new count
even in these cases, where other computational techniques are available.
\end{abstract}

\maketitle

\section{Introduction: the numbers $N^{\Delta,\delta}$}
%\subsection{The set-up of the problem in $\tor$}
Let $\Delta\subset\R^2$ be a convex polygon with integer vertices.
It defines a finite-dimensional linear system $\pp\LLL$ of curves in
$\tor$,
where $\C^*=\C\setminus\{0\}$.
These curves are the zero loci in $\tor$ of the (Laurent) polynomials
$$f(z,w)=\sum\limits_{(j,k)\in\Delta\cap\Z^2}a_{jk}z^jw^k,$$
$a_{jk}\in\C$.
The polynomials $f$ themselves form the vector space $\LLL$.
Recall that {\em the Newton polygon} of $f$ is
$\operatorname{Convex hull}\{(j,k)\ |\ a_{jk}\neq 0\}$.
Thus $\LLL$ contains polynomials
whose Newton polygon is contained in $\Delta$.
Clearly $\pp\LLL$ is a complex projective space of dimension
$$m=\#(\Delta\cap\Z^2)-1.$$ Curves with the Newton polygon $\Delta$
form an open dense set $\U\subset\pp\LLL$.
A generic curve in $\LLL$ is a smooth curve of genus
$$l=\#(\Int\Delta\cap\Z^2).$$

Let $C\in\pp\LLL$ be a curve. Even if $C$ is not irreducible
we can define its genus $g(C)\in\Z$.
Consider the decomposition
$C=C_1\cup\dots\cup C_n$ into the irreducible components $C_j$.
We define $g(C)=\sum\limits_{j=1}^n g(C_j)+1-n$.
%This definition agrees with the Euler characteristic formula.
Note that this definition of genus allows for negative values
(cf. \cite{CH}).
If $C$ is singular then its genus is strictly smaller than $l$.

The curves of genus $l-\delta$ and with the Newton polygon $\Delta$
form a subvariety $\Sigma^\circ_\delta\subset\U$
of dimension $m-\delta$.
Let $\Sigma_\delta\subset\LLL$ be the projective closure of
$\Sigma^\circ_\delta$.
We define $N^{\Delta,\delta}$ to be the degree
%(i.e. the intersection
%number with a projective subspace of codimension $m-\delta$)
of the $(m-\delta)$-dimensional subvariety $\Sigma_\delta$ in $\LLL$.
The degree is the intersection number with
a projective subspace of codimension $m-\delta$.
Curves from $\LLL$ passing through a point $z\in\tor$ form
a hyperplane.

The number $N^{\Delta,\delta}$
has the following enumerative interpretation.
Let $z_1,\dots,z_{m-\delta}\in\tor$ be generic points.
The number $N^{\Delta,\delta}$ equals to the number of
algebraic curves of Newton polygon $\Delta$ and genus $l-\delta$
passing through $z_1,\dots,z_{m-\delta}$.
%Note that this dimension is not always pure.
%Let $\Delta'\subset\Delta$ be a subpolygon disjoint from a side of
%$\Delta$.
Note that $N^{\Delta,0}=1$ for any $\Delta$. These numbers get
more interesting when $\delta>0$.

%\subsection{Compactification:
%the set-up of the same problem in $\C T_\Delta$}
\begin{rmk}
\label{comp}
Another way to look at the same problem is to consider
the compactification of the torus $\tor$.
Recall that
the polygon $\Delta$ defines a compact toric surface $\C T_\Delta$,
see e.g. \cite{GKZ}.
(Some readers may be more familiar with the definition of toric
surfaces by fans, in our case the fan is formed by the dual
cones at the vertices of $\Delta$.)
The surface $\C T_\Delta$ may have isolated singularities
that correspond to some vertices of $\Delta$.

In addition to a complex structure (which depends only on the
dual fan) the polygon $\Delta$ defines a holomorphic linear bundle
$\HH$ over $\C T_\Delta$.
We have a canonical identification $\Gamma(\HH)=\LLL$,
where $\Gamma(\HH)$ is the space of the sections of $\HH$.
The projective space $\pp\LLL$ can be also considered
as the space of all holomorphic curves in $\C T_\Delta$
such that their homology class is Poincar\'e dual to
$c_1(\HH)$.
%
%Let $\U'\subset\p\LLL$ be formed by the curves that do not
%have components contained in $\C T_\Delta\setminus\tor$.
%Note that $\U'\subset U$.
%If $\bar{C}\subset\C T_\Delta$ is a curve
%and $\tilde{C}\to \bar{C}$ is its normalization then we define
%$g(\bar{C})=\frac12 (2-\chi(\tilde{C})).$ Note that if $\bar{C}\in\U'$
%then $g(\bar{C})=g(C)$, where
%$C=\bar{C}\cap\tor$ (its genus is defined in the previous subsection).
%
%
%The open dense subset $\U\subset\p\LLL$ is formed by those
%curves in $\C T_\Delta$ that do not pass through the vertices
%of the toric surface $\C T_\Delta$. In particular, such curves
%do not pass through the singular points of $\C T_\Delta$.
%
%Let $z_1,\dots,z_{m_\delta}\in \C T_{\Delta}$ be generic points.
The number $N^{\Delta,\delta}$ is the number of holomorphic
curves $\bar{C}\subset\C T_\Delta$ such that
$z_1,\dots,z_{m_\delta}\in\bar{C}$,
the homology class $[\bar{C}]$ is dual to $c_1(\HH)$,
the Euler characteristic of the normalization of $\bar{C}$
is $2-2(l-\delta)$ and no irreducible component of
$\bar{C}$ is contained in
$\C T_\Delta\setminus\tor$.
\end{rmk}

\begin{rmk}
Note that in the set-up of Remark \ref{comp}
the number $N^{\Delta,\delta}$ appears
related to the Gromov-Witten invariant of $\C T_\Delta$ (see \cite{KM})
corresponding to $c_1(\HH)$.
The difference is that the corresponding Gromov-Witten
invariant also has a contribution from the curves of genus
$l-\delta$ which have components contained in $\C T_\Delta\setminus\tor$.
This contribution is zero (by the dimension reasons)
if $\C T_\Delta$ is smooth and does not have exceptional divisors.
Thus if $\Delta=\Delta_d=\operatorname{Convex Hull}\{(0,0),(d,0),
(0,d)\}$ or $\Delta=[0,r]\times[0,s]$ then the number
$N^{\Delta,\delta}$ is the multicomponent
Gromov-Witten invariant of genus $(l-\delta)$ and degree $d$ in $\cp^2$
or of bidegree $(r,s)$ in $\cp^1\times\cp^1$.
\end{rmk}

\begin{exam}
Suppose $\Delta=\Delta_d$ so that $\C T_\Delta=\cp^2$.
%the most famous setup for the enumerative problem.
We have $m=\frac{d(d+3)}{2}$ and $l=\frac{(d-1)(d-2)}{2}$.
The number $N^{\Delta_d,\delta}=N_{g,d}$ is the number of
genus $g$, degree $d$ (not necessarily irreducible) curves
passing through $3d-1+g$ generic points in
$\cp^2$, $g=\frac{(d-1)(d-2)}{2}-\delta$.

%Some first non-trivial numbers in this series are $N_{0,3}=12$,
%$N_{2,4}=27$, $N_{1,4}=225$, $N_{0,4}=675$, $N_{-1,2}=3$,
%$N_{-1,3}=21$.
%
The formula $N^{\Delta_d,1}=3(d-1)^2$
is well-known as the degree of the discriminant (cf. \cite{GKZ}).
An elegant recursive formula for the number of irreducible
rational curves (the one-component part of $N_{0,d}$)
was found by Kontsevich \cite{KM}.
An algorithm for computing $N_{g,d}$ for arbitrary $g$
is due to Caporaso and Harris \cite{CH}.
See \cite{V} for computations for some other rational surfaces,
in particular, the Hirzebruch surfaces
(this corresponds to the case when $\Delta$ is a trapezoid).
\end{exam}

\section{Lattice paths and their multiplicities} \label{paths}
A path $\gamma:[0,n]\to\R^2$, $n\in\N$, is called a {\em lattice path}
if $\gamma|_{[j-1,j]}$, $j=1,\dots,n$ is an affine-linear map
and $\gamma(j)\in\Z^2$, $j\in 0,\dots,n$.
Clearly, a lattice path is determined by its values at
the integer points.
%To any such path we associate a sequence of its normal vectors
%$\gamma_j$, $j=1,\dots,n$. Each $\gamma_j$ is obtained by
%rotating the vector connecting $\gamma(j-1)$ and $\gamma(j)$
%by 90 degrees counterclockwise. Note that each $\gamma_j$ is an
%integer vector. Note also that the point
%$\gamma(0)$ and the vectors $\gamma_j$ determine
%the lattice path $\gamma$.
%
Let us choose an auxiliary linear map $\lambda:\R^2\to\R$ that is irrational,
i.e. such that $\lambda|_{\Z^2}$ is injective.
%The map $\alpha$ defines an order on the lattice points of $\Delta$.
Let $p,q\in\Delta$ be the vertices where $\alpha|_\Delta$ reaches
its minimum and maximum respectively.
A lattice path is called {\em $\lambda$-increasing}
if $\lambda\circ\gamma$ is increasing.

%We say that a vector $\beta$ is $\lambda$-positive
%\begin{defn}
%We define the {\em $\lambda$-multiplicity} of
%a sequence $\beta_1,\dots,\beta_k\in\Z^2$ inductively
%by $k$.
%\begin{itemize}
%\item If $k=1$ then
%\end{itemize}
%\end{defn}

The points $p$ and $q$ divide the boundary $\dd\Delta$ into
two increasing lattice paths
$$\alpha^+:[0,n_+]\to\dd\Delta\ \  \text{and}\ \
\alpha^-:[0,n_-]\to\dd\Delta.$$
We have $\alpha_+(0)=\alpha_-(0)=p$,
$\alpha_+(n_+)=\alpha_-(n_-)=q$, $n_++n_-=m-l+3$.
To fix a convention we assume that $\alpha_+$ goes clockwise
around $\dd\Delta$ wile $\alpha_-$ goes counterclockwise.
%This gives us the vectors $\alpha^+_1,\dots,\alpha^+_{n_+},
%\alpha^-_1,\dots,\alpha^-_{n_-}\in\Z^2$.

Let
$\gamma:[0,n]\to\Delta\subset\R^2$
be an increasing lattice path such that $\gamma(0)=p$
and $\gamma(n)=q$.
The path $\gamma$ divides $\Delta$ into two closed regions:
$\Delta_+$ enclosed by $\gamma$ and $\alpha_+$
and $\Delta_-$ enclosed by $\gamma$ and $\alpha_-$.
Note that the interiors of $\Delta_+$ and $\Delta_-$
do not have to be connected.

We define the positive (resp. negative)
multiplicity $\mu_\pm(\gamma)$
of the path $\gamma$ inductively.
We set $\mu_\pm(\alpha_\pm)=1$.
If $\gamma\neq\alpha_\pm$
%then $\Delta_\pm$ is a polygon
%with a non-empty interior.
%Let
then we take $1\le k\le n-1$ to be the smallest number
such that $\gamma(k)$ is a vertex of $\Delta_\pm$
with the angle less than $\pi$ (so that $\Delta_\pm$
is locally convex at $\gamma(k)$).

If such $k$ does not exist we set $\mu_\pm(\gamma)=0$.
If $k$ exist we consider two other increasing lattice paths connecting
$p$ and $q$
$\gamma':[0,n-1]\to\Delta$ and $\gamma'':[0,n]\to\R^2$.
We define $\gamma'$ by $\gamma'(j)=\gamma(j)$ if $j<k$ and
$\gamma'(j)=\gamma(j+1)$ if $j\ge k$.
We define $\gamma''$ by $\gamma''(j)=\gamma(j)$ if $j\neq k$ and
$\gamma''(k)=\gamma(k-1)+\gamma(k+1)-\gamma(k)\in\Z^2$.
We set
$$\mu_\pm(\gamma)=2\operatorname{Area}(T)\mu_\pm(\gamma')+
\mu_\pm(\gamma''),$$
where $T$ is the triangle with the vertices $\gamma(k-1)$,
$\gamma(k)$ and $\gamma(k+1)$.
The multiplicity is always integer since the area of a lattice
triangle is half-integer.

Note that it may happen that $\gamma''(k)\notin\Delta$.
In such case we use a convention  $\mu_\pm(\gamma'')=0$.
We may assume that $\mu_\pm(\gamma')$ and $\mu_\pm(\gamma'')$
is already defined since the area of $\Delta_\pm$ is smaller
for the new paths. Note that $\mu_\pm=0$ if $n<n_\pm$
as the paths $\gamma'$ and $\gamma''$ are not longer than $\gamma$.

We define {\em the multiplicity of the path $\gamma$}
as the product $\mu_+(\gamma)\mu_-(\gamma)$.
Note that the multiplicity of a path connecting two
vertices of $\Delta$ does not depend on $\lambda$.
We only need $\lambda$ to determine whether a path
is increasing.

\begin{exa}
Consider the path $\gamma:[0,8]\to\Delta_3$ depicted
on the extreme left of Figure \ref{mult+}.
This path is increasing
with respect to $\lambda(x,y)=x-\epsilon y$, where
$\epsilon>0$ is very small.

\begin{figure}[h]
\centerline{\psfig{figure=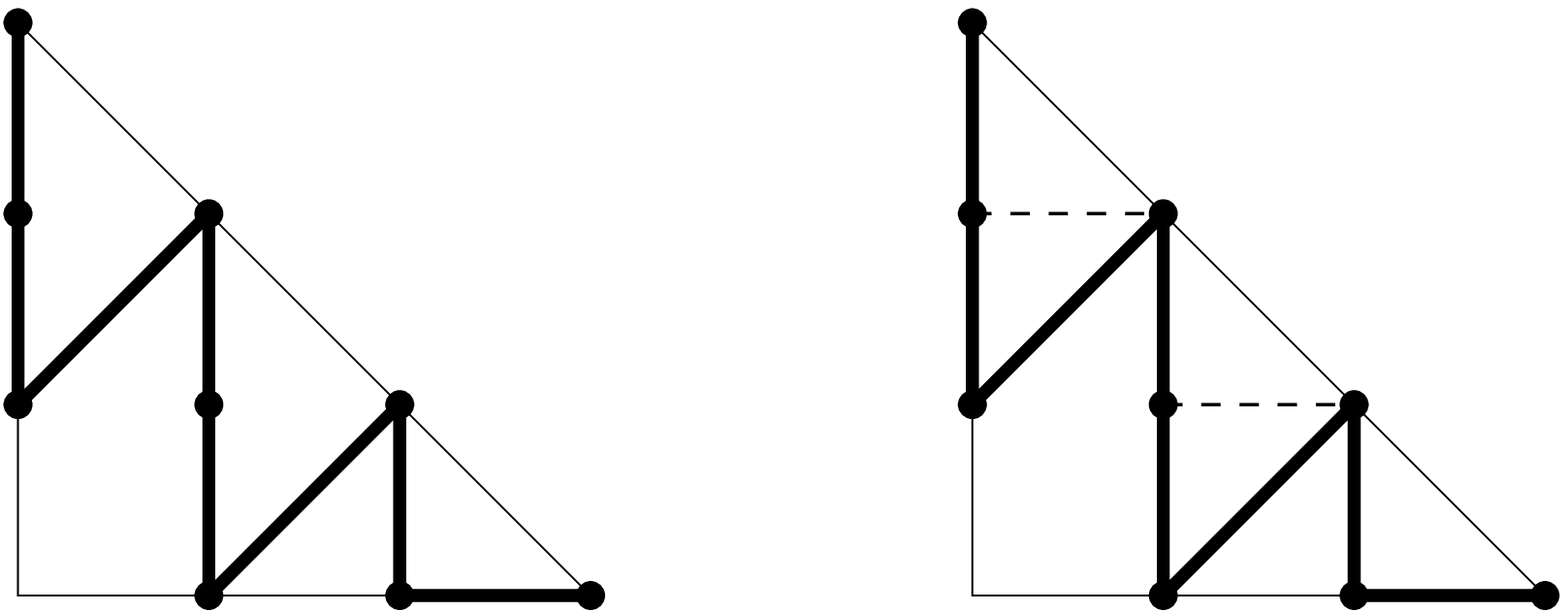,height=0.85in,width=2.2in}
\hspace{0.5in}
\psfig{figure=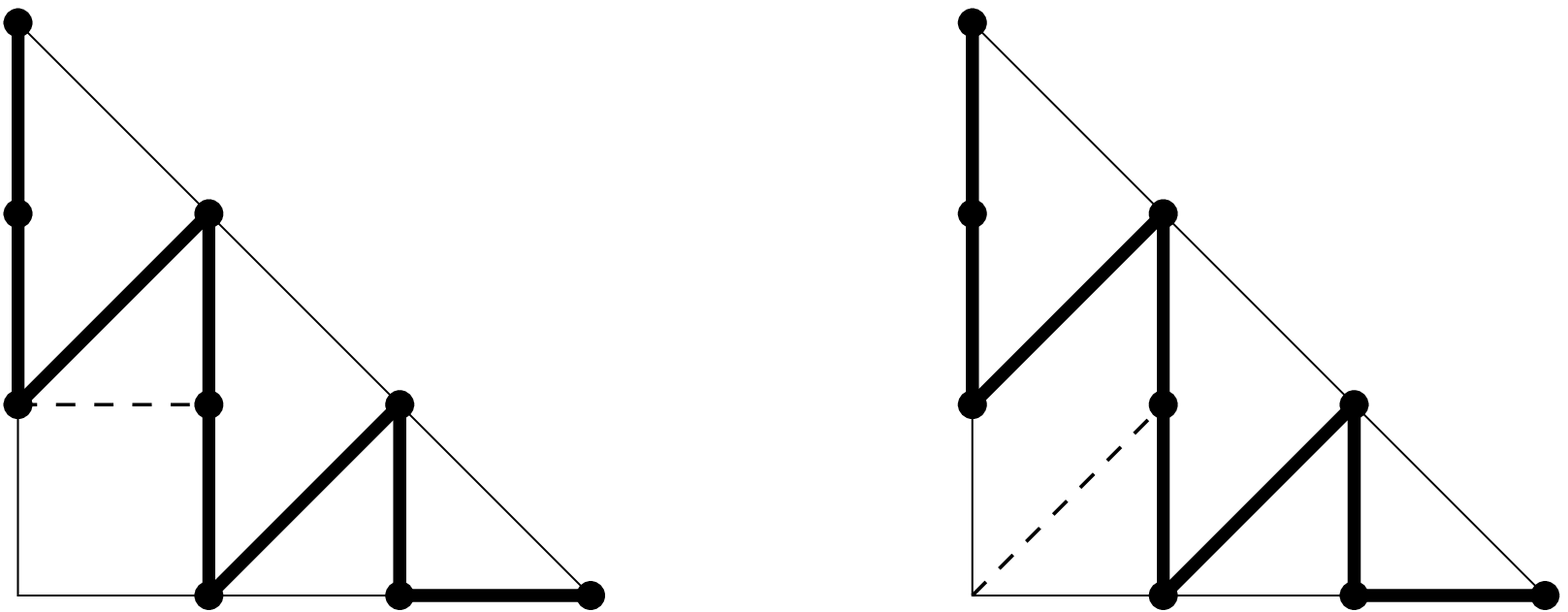,height=0.85in,width=2.2in}}
\caption{\label{mult+} A path $\gamma$ with $\mu_+(\gamma)=1$
and $\mu_-(\gamma)=2$.}
\end{figure}

Let us compute $\mu_+(\gamma)$. We have $k=2$ as $\gamma(2)=(0,1)$
is a locally convex vertex of $\Delta_+$.
We have $\gamma''(2)=(1,3)\notin\Delta_3$ and thus
$\mu_+(\gamma)=\mu_+(\gamma')$, since $\operatorname{Area}(T)=\frac12$.
Proceeding further we get
$\mu_+(\gamma)=\mu_+(\gamma')=\dots=\mu_+(\alpha_+)=1$.
%The second picture of Figure \ref{mult+} shows the subdivision

%\begin{figure}[h]
%\centerline{\psfig{figure=mult-.eps,height=1in,width=2.7in}}
%\caption{\label{mult-}  $\mu_-(\gamma)=2$, since
%$\mu_-(\gamma')=1$ and
%$\mu_-(\gamma'')=1$.}
%\end{figure}

Let us compute $\mu_-(\gamma)$. We have $k=3$ as $\gamma(3)=(1,2)$
is a locally convex vertex of $\Delta_-$.
We have $\gamma''(3)=(0,0)$ and $\mu_-(\gamma'')=1$.
To compute $\mu_-(\gamma')=1$ we note that $\mu_-((\gamma')')=0$
and $\mu_-((\gamma')'')=1$.
Thus the full multiplicity of $\gamma$ is 2.
\end{exa}

\section{The formula}
In the previous section we fixed
an auxiliary linear function $\lambda:\R^2\to\R$
which
determines the extremal vertices $p,q$ of $\Delta$.

\begin{thm}
\label{thm1}
The number $N^{\Delta,\delta}$ equals to the
number (counted with multiplicities)
of $\lambda$-increasing lattice paths
$[0,m-\delta]\to\Delta$ connecting $p$ and $q$.
\end{thm}

This theorem is proved in \cite{M} (to appear).
The proof is based on the application of the so-called
{\em tropical algebraic geometry} (see e.g. Chapter 9 of \cite{S}).
The relation between the classical enumerative problem and
the corresponding tropical problem is provided by
passing to the ``large complex limit" as suggested
by Kontsevich (see \cite{KS} for these ideas in a more general setting).

Note that an immediate corollary of Theorem \ref{thm1} is that
the number of $\lambda$-increasing lattice paths
of a fixed length does not depend on the choice of $\lambda$.

\begin{exa}\label{ecusp}
Let us compute $N^{\Delta,1}=5$ for the polygon $\Delta$
depicted on Figure \ref{cusp} in two different ways.
Using $\lambda(x,y)=-x+\epsilon y$ for a small $\epsilon>0$
we get the left two paths depicted on Figure \ref{cusp}.
Using $\lambda(x,y)=x+\epsilon y$ we get the three right paths.
The corresponding multiplicities are shown under the path.
All other $\lambda$-increasing paths have zero multiplicity.
\begin{figure}[h]
\centerline{\psfig{figure=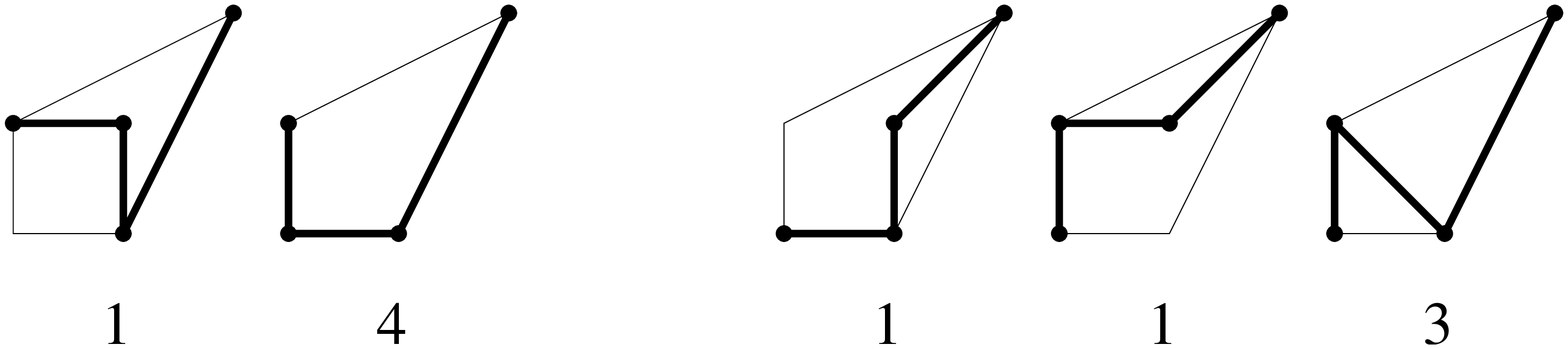,height=0.85in,width=3.9in}}
\caption{\label{cusp} Computing $N^{\Delta,1}=5$ in two different ways.}
\end{figure}
\end{exa}

In the next two examples we use $\lambda(x,y)=x-\epsilon y$
as the auxiliary linear function.

\begin{exa}\label{edeg3}
Figure \ref{deg3} shows a computation of
the well-known number $N^{\Delta_3,1}=N_{0,3}$.
This is the number of rational cubic curves through 8 generic
points in $\cp^2$.
\begin{figure}[h]
\centerline{\psfig{figure=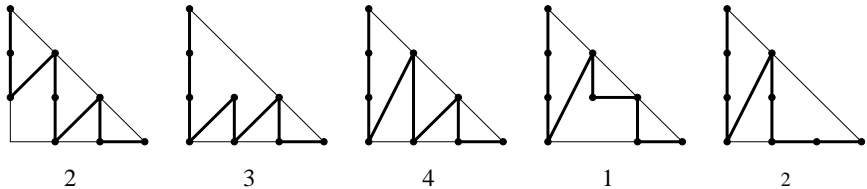,height=0.95in,width=4.5in}}
%\hspace{0.1in}
%\psfig{figure=deg3alt.eps,height=0.6in,width=2.4in}}
\caption{\label{deg3} Computing $N_{0,3}=12$.}
\end{figure}
%\begin{figure}[h]
%\centerline{\psfig{figure=deg3alt.eps,height=1.2in,width=5in}}
%\caption{\label{deg3alt} A different computation of $N_{0,3}=12.$}
%\end{figure}
\end{exa}
%In the previous example we had the same set of multiplicities
%of the $\lambda$-increasing paths for different choices of $\lambda$,
%12=4+3+2+2+1. This is not true in general.
%E.g. if $\Delta=[0,3]\times [0,2]$ (i.e. we count curves
%of bidegree $(3,2)$ in $\p^1\times\p^1$) then the choice
%$\lambda(x,y)=x-\epsilon y$ yields the computation
%$N^{\Delta,1}=4+4+2+2+2+2+2+2$ while the choice
%$\lambda(x,y)=y-\epsilon x$ yields $N^{\Delta,1}=4+4+3+3+3+3$.
%
\begin{exa}\label{edeg4}
Figure \ref{deg4} shows a computation
of a less well-known number $N^{\Delta_4,2}=N_{1,4}$.
This is the number of genus 1 quartic curves through
12 generic points in $\cp^2$.
\begin{figure}[h]
\centerline{\psfig{figure=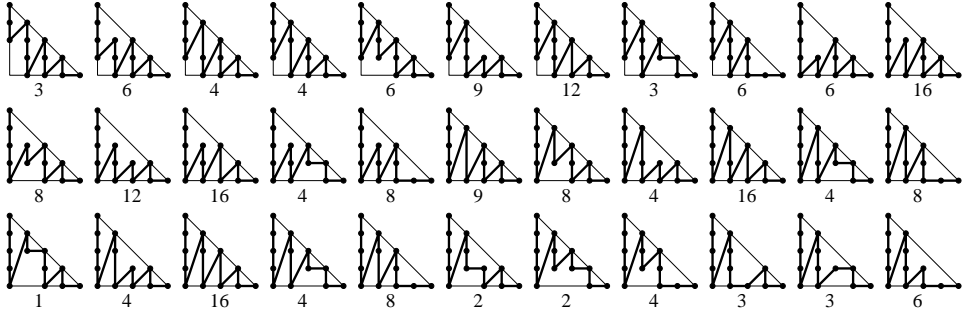,height=1.6in,width=5in}}
\caption{\label{deg4} Computing $N_{1,4}=225$.}
\end{figure}
\end{exa}

\section{Real aspects of the count}
Suppose that $z_1,\dots,z_{m-\delta}\in\rtor\subset\tor$
are generic real points.
We may ask how many of the $N^{\Delta,\delta}$
relevant complex curves are real, i.e. defined over $\R$.
Note that this number depends on the configuration of real points.

Theorem \ref{thm1} can be modified to give the relevant count
of real curves. In order to do this we need to
define the real multiplicity of a lattice path
$\gamma:[0,n]\to\Delta$ connecting the vertices $p$ and $q$.
%and equipped with signs $\sigma_j\in\Z_2\oplus\Z_2$, $j=1,\dots,n$.
We introduce the sequence
of the pairs of signs $\sigma_1,\dots,\sigma_{n}\in\Z_2\oplus\Z_2$
(this sequence will record the quadrants of the points $z_j\in\rtor$).
%where $\sigma_j\in\Z_2\oplus\Z_2$ given by the quadrants of
%the points $z_j$.
The sign $\sigma_j$ is prescribed to the edge $\gamma[j-1,j]$.
We make a convention that $\sigma_j$ and $\sigma'_j$
are equivalent if
$\sigma_j-\sigma'_j\equiv \gamma(j)-\gamma(j-1)\pmod{2}.$
%(Thus to an "odd" edge of $\gamma$ we associate a pair of equivalent
%sign choices while for an "even" edge the sign choice is unique even
%up to the equivalence.)

We set
\begin{equation}
\label{remult}
\mu^{\R}_\pm(\gamma)=a(T)\mu^{\R}_\pm(\gamma')+
\mu^{\R}_\pm(\gamma'').
\end{equation}
The definition of the new paths $\gamma'$, $\gamma''$ and the
triangle $T$ is the same as in section \ref{paths}. The sign
sequence for $\gamma''$ is $\sigma''_j=\sigma_j,j\neq k,k+1$,
$\sigma''_k=\sigma_{k+1}$, $\sigma''_{k+1}=\sigma_k$. The sign
sequence for $\gamma'$ is $\sigma'_j=\sigma_j,j< k$,
$\sigma'_j=\sigma_{j+1}, j>k$.
%while $\sigma'_{k}$ is defined below.
We define the sign $\sigma'_{k}$ and the function $a(T)$ as
follows.
%according to the following three cases.

%\begin{description}
%\item[If all sides of $T$ are odd]
%Suppose that all three vectors (the sides of $T$)
%$\gamma(k+1)-\gamma(k-1)$,
%$\gamma(k-1)-\gamma(k)$ and $\gamma(k+1)-\gamma(k)$
%are not divisible by 2 in $\Z^2$.
%%Suppose that $\operatorname{Area}(t)$ is not integer
%%(i.e. all three sides of $T$ have "odd" length).
\begin{itemize}
\item If all sides of $T$ are odd we set $a(T)=1$ and define the
sign $\sigma'_k$ (up to the equivalence) by the condition that the
three equivalence classes of $\sigma_k$, $\sigma_{k+1}$ and
$\sigma'_k$ do not share a common element.
%\item[If all sides of $T$ are even]
%Suppose that all three vectors
%$\gamma(k+1)-\gamma(k-1)$,
%$\gamma(k-1)-\gamma(k)$ and $\gamma(k+1)-\gamma(k)$
%are divisible by 2.
%%(i.e. all three sides of $T$ have "even" length).
\item If all sides of $T$ are even we set $a(T)=0$ if
$\sigma_{k-1}\neq\sigma_k$. In this case we can ignore $\gamma'$
(and its sequence of signs). We set $a(T)=4$ if
$\sigma_{k}=\sigma_{k+1}$. In this case we define
$\sigma'_k=\sigma_k=\sigma_{k+1}$. \item Otherwise
%Suppose that none of the above two conditions hold
%(then exactly one of the three vectors is divisible by 2).
we set $a(T)=0$ if the equivalence classes of $\sigma_k$ and
$\sigma_{k+1}$ do not have a common element. We set $a(T)=2$ if
they do. In the latter case we define the equivalence class of
$\sigma'_k$ by the condition that $\sigma_k$, $\sigma_{k+1}$ and
$\sigma'_k$ have a common element. There is one exception to this
rule. If the even side is
%the odd vectors are $\gamma(k-1)-\gamma(k)$ and
%$\gamma(k+1)-\gamma(k)$ and the even vector is
$\gamma(k+1)-\gamma(k-1)$ then there are two choices for
$\sigma'_k$ satisfying the above condition. In this case we
replace $a(T)\mu^{\R}_\pm(\gamma')$ in \eqref{remult} by the sum
of the two multiplicities of $\gamma'$ equipped with the two
allowable choices for $\sigma'_k$ (note that this agrees with
$a(T)=2$ in this case).
%\end{description}
\end{itemize}
Similar to section \ref{paths} we define
$\mu^{\R}_\pm(\alpha_\pm)=1$ and
$\mu^{\R}(\gamma)=\mu^{\R}_+(\gamma)\mu^{\R}_-(\gamma)$. As before
$\lambda:\R^2\to\R$ is a linear map injective on $\Z^2$ and $p$
and $q$ are the extrema of $\lambda|_\Delta$.
\begin{thm} \label{thm2}
For any choice of $\lambda$ and $\sigma_j, j=1,\dots,m-\delta$
there exists a configuration of $m-\delta$ of generic points in
the respective quadrants
%$(\R_{>0})^2\subset\rtor$
such that the number of real curves among the $N^{\Delta,\delta}$
relevant complex curves is equal to the number of
$\lambda$-increasing lattice paths $\gamma:[0,m-\delta]\to\Delta$
connecting $p$ and $q$ counted with multiplicities $\mu^{\R}$.
\end{thm}
%We use the notation $\R_+=\{x>0\}$ so that $(\R_+)^2$ stands
%for the positive quadrant of the torus $\rtor$.

%In the examples below we use the choice $\sigma_j=(+,+)$ so all
%the points $z_j$ are in the positive quadrant
%$(\R_{>0})^2\subset\rtor$.

\begin{exa} Here we use the choice $\sigma_j=(+,+)$ so all
the points $z_j$ are in the positive quadrant
$(\R_{>0})^2\subset\rtor$. The first count of $N^{\Delta,1}$ from
Example \ref{ecusp} gives a configuration of 3 real points with 5
real curves. The second count gives a configuration with 3 real
curves as the real multiplicity of the last path is 1. Note also
that the second path on Figure \ref{cusp} changes its real
multiplicity if we reverse its direction.

Example \ref{edeg3} gives a configuration of 9 generic points in
$\rp^2$ with all 12 nodal cubics through them real. Example
\ref{edeg4} gives a configuration of 12 generic points in $\rp^2$
with 217 out of the 225 quartics of genus 1 real. The path in the
middle of Figure \ref{deg4} has multiplicity 9 but real
multiplicity 1. A similar computation shows that there exists a
configuration of 11 generic points in $\rp^2$ such that 564 out of
the 620 irreducible quartic through them are real.
\end{exa}

\begin{rmk}
\label{rgw}
%One can associate another real multiplicity to a path $\gamma$ by
%replacing \eqref{remult} with
Real nodal curves have three types of nodes: hyperbolic, elliptic
and imaginary. Theorem \ref{thm2} can be refined to count curves
with different types of nodes separately. In accordance with
\cite{W} let us prescribe a sign $(-1)^e$ to a real nodal curve,
where $e$ is the number of its elliptic nodes.
%One can modify Theorem \ref{thm2} to compute the
%corresponding algebraic number
%of real curves via a configuration of generic real points by replacing
%$\mu^{\R}$ with another multiplicity $\nu^{\R}$.
%a real curve is counted with the sign equal to $(-1)^e$
To compute the corresponding algebraic number of curves we
introduce the multiplicity $\nu^{\R}$ by replacing \eqref{remult}
with $\nu^{\R}_\pm(\gamma)=b(T)\nu^{\R}_\pm(\gamma')+
\nu^{\R}_\pm(\gamma'').$ Here we define $b(T)=0$ if at least one
side of $T$ is even and $b(T)=(-1)^{\#(\Int T\cap\Z^2)}$
otherwise. It can be shown with the help of this formula and a
combinatorial observation made by Itenberg, Kharlamov and Shustin
(to appear)
%$\lambda(x,y)=x-\epsilon y$ and
%$\Delta=\Delta_d$
that in the case $\Delta=\Delta_d$ the algebraic number of {\em
irreducible} curves counted by $\nu^{\R}$ is positive for any
genus $0\le g\le\frac{(d-1)(d-2)}{2}$ if $\lambda(x,y)=y-\epsilon
x$.

Note that unlike $\mu^{\R}$ the multiplicity $\nu^{\R}$ does not
depend on the quadrant choices $\sigma_j$. Furthermore, in
\cite{W} Welschinger stated that if $g=0$ then this number is
independent of the configuration of generic real points. Corollary
1.2 of \cite{W} combined with Remark \ref{rgw} implies the
following statement (which answers a question asked e.g. by
Rokhlin and Kharlamov). {\em For any configuration of generic
$3d-1$ points in $\rp^2$ there exists a real rational curve of
degree $d$ passing through this configuration.}
\end{rmk}

\end{document}